\documentclass[11pt,a4paper]{amsart}
\usepackage[utf8]{inputenc}
\usepackage[english]{babel}
\usepackage[dvipsnames]{xcolor}
\usepackage[all]{xy}
\usepackage{hyperref}
\usepackage{tikz}
\usepackage{tikz-cd}
\usepackage{soul,color}
\usepackage{comment}
\usepackage[colorinlistoftodos]{todonotes}
\usepackage[width=\textwidth,font=footnotesize,labelfont=bf, aboveskip=30pt]{caption}
\usepackage[symbol]{footmisc}
\usepackage[OT2,T1]{fontenc}
\usepackage{amsthm}
\DeclareSymbolFont{cyrletters}{OT2}{wncyr}{m}{n}
\DeclareMathSymbol{\Sha}{\mathalpha}{cyrletters}{"58}

\newcommand\sti[1]{\todo[inline,color=green]{#1}}
\newcommand\vni[1]{\todo[inline,color=cyan]{#1}}
\newcommand\tksi[1]{\todo[inline,color=yellow]{#1}}

\newcommand*{\fullref}[1]{\hyperref[{#1}]{\ref*{#1}. \nameref*{#1}}} 
\theoremstyle{remark}

\newcommand{\cH}{\mathrm{H}}

\usepackage[osf]{mathpazo}
\usepackage{charter}
\usepackage[a4paper, top=3.2cm, bottom=3.2cm,left=2.8cm, right=2.8cm, heightrounded,bindingoffset=0mm]{geometry}

\usepackage{fancyhdr}

\usepackage{amsfonts}
\usepackage{amssymb}
\usepackage{extarrows}
\usepackage{cleveref}
\usepackage{amsthm}
\usepackage{amscd}
\usepackage{amsmath}
\usepackage{mathtools}
\usepackage{mathrsfs}
\usepackage{float}

\usepackage{color}	
\definecolor{due}{RGB}{0,76,147}

\usepackage{graphicx}	
\usepackage{multicol}	
\usepackage{wrapfig}
\usepackage[english]{babel} 
\usepackage[utf8]{inputenc}
\usepackage[T1]{fontenc}
\usepackage{enumitem}

\usepackage{longtable}
\usepackage{cite}

\theoremstyle{definition}
\newtheorem{defi}{Definition}[section]
\newtheorem{setting}[defi]{Setting}
\theoremstyle{plain}

\newtheorem*{thm*}{Theorem}

\newtheorem*{mainthm}{Main Theorem}

\newtheorem{prop}[defi]{Proposition}

\newtheorem{lemma}[defi]{Lemma}
\theoremstyle{remark}

\newtheorem{rmk}[defi]{Remark}
\newtheorem{ex}[defi]{Example}

\theoremstyle{definition}

\newcommand{\kl}{{\mathcal L}}
\newcommand{\km}{{\mathcal M}}

\newcommand{\kx}{{\mathcal X}}

\newcommand{\IG}{{\mathbb G}}

\newcommand{\Z}{{\mathbb Z}}

\newcommand{\hx}{\hat{x}}
\newcommand{\hy}{\hat{y}}
\newcommand{\hz}{\hat{z}}
\newcommand{\hw}{\hat{w}}
\newcommand{\AutcX}{\mathcal{A}ut_{\mathcal{X}}}

\newcommand{\DXa}{\operatorname{D}^b(X,\alpha)}

\newcommand{\Aut}{\operatorname{Aut}}

\newcommand{\Pic}{\operatorname{Pic}}

\newcommand{\ol}{\overline}

  \makeatletter
\newcommand{\xdashrightarrow}[2][]{\ext@arrow 0359\rightarrowfill@@{#1}{#2}}
\newcommand{\xdashleftarrow}[2][]{\ext@arrow 3095\leftarrowfill@@{#1}{#2}}
\newcommand{\xdashleftrightarrow}[2][]{\ext@arrow 3359\leftrightarrowfill@@{#1}{#2}}
\def\rightarrowfill@@{\arrowfill@@\relax\relbar\rightarrow}
\def\leftarrowfill@@{\arrowfill@@\leftarrow\relbar\relax}
\def\leftrightarrowfill@@{\arrowfill@@\leftarrow\relbar\rightarrow}
\def\arrowfill@@#1#2#3#4{%
  $\m@th\thickmuskip0mu\medmuskip\thickmuskip\thinmuskip\thickmuskip
   \relax#4#1
   \xleaders\hbox{$#4#2$}\hfill
   #3$%
}
\makeatother

\numberwithin{equation}{section}
\usepackage{hyperref}
\hypersetup{bookmarksnumbered=true, linkcolor=due, citecolor=Green, urlcolor=black,}
\usepackage{amssymb}
\newcommand{\bk}{\Bbbk}
\usepackage{mathrsfs}
\providecommand{\kx}{}\renewcommand{\kx}{\mathscr{X}}
\providecommand{\AutcX}{}\renewcommand{\AutcX}{\mathscr{A}ut_{\kx}}
\makeatletter
\@ifundefined{rmk}{%
  \theoremstyle{remark}%
  \@ifundefined{c@prop}{\newtheorem{rmk}{Remark}[section]}{\newtheorem{rmk}[prop]{Remark}}%
}{}
\makeatother

\begin{document}
	\title[Twisted Genus 1 curves]{Twisted Fourier--Mukai partners of genus 1 curves}
	\author{Tanya Kaushal Srivastava}
	\address{Department of Mathematics\\
    Indian Institute of Technology Gandhinagar, Palaj Campus, Gandhinagar, Gujarat, India}
\email{\url{ks.tanya@iitgn.ac.in}}
\author{Sofia Tirabassi}
\address{Department of Mathematics\\ Stockholm University, Kraftriket hus 6, Stockholm, Sweden}
\email{\url{tirabassi@math.su.se}}
\author{Naïa Vivares}
\address{Institut Math\'ematiques de Jussieu-Paris Rive Gauche \\ Universit\'e Paris Cit\'e , 8 Place Aur\'elie N\'emours, Paris, France}
\email{\url{vivares@imj-prg.fr}}

\begin{abstract}
We study twisted derived equivalence between genus 1 curves over an arbitrarly perfect field $\bk$.
\end{abstract}
\clearpage\maketitle 

\setcounter{tocdepth}{1}
\tableofcontents
\section{Introduction}
In this paper we study twisted Fourier--Mukai partners of genus 1 curves over a perfect field $\bk$.\par
We recall that a \emph{twisted variety} over $\bk$ is given by a pair $(X,\alpha)$, with $X$  a smooth projective geometrically connected variety over $\bk$ and $\alpha\in \operatorname{Br}(X)$ an element of its Brauer group. Given a twisted variety, there is a $\IG_m$-gerbe $\pi\colon\kx\to X$ with class $\alpha$, unique up to isomorphism. We write $\DXa=\operatorname{D}^b(\kx)^{(1)}$ for the derived category of $1$-twisted sheaves.  We say that two varieties $X$ and $Y$ are \emph{twisted derived equivalent} if there are two Brauer classes $\alpha \in \operatorname{Br}(X)$ and $\beta\in\operatorname{Br}(Y)$ and  an exact equivalence 
\[\Phi:\DXa\rightarrow \operatorname{D}^b(Y,\beta).\]
We say that two varieties $X$ and $Y$ are \emph{derived isogenous} if there are a finite sequence of varieties
$$
X\simeq X_0,\; X_1,\ldots ,X_n\simeq Y
$$
and Brauer classes $\beta_i\in\operatorname{Br}(X_i)$ and $\alpha_{i+1}\in\operatorname{Br}(X_{i+1})$ such that $\operatorname{D}^b(X_i,\beta_i)\simeq \operatorname{D}^b(X_{i+1},\alpha_{i+1})$ for $0\le i<n$.
Investigating the hidden relationships between twisted derived equivalent or derived isogenous varieties is a question that has raised, both recently and less recently, much interest (see for example \cite{Add18, li2025derived, li2025point, lane2026twisted}).\par
In this paper we completely characterize derived isogenous genus 1 curves, generalizing the results of \cite{antieau2017derived} to the twisted setting.

\begin{mainthm}
Let $\bk$ be a perfect field and let $C,C'$ be genus 1 curves over $\bk$. The following are equivalent:
\begin{enumerate}
\item $C$ and $C'$ are derived isogenous;
\item the Jacobians $J(C)$ and $J(C')$ are isomorphic over $\bk$;
\item $C$ and $C'$ are connected by a chain (of length two) of twisted derived equivalences passing through $J(C)$.
\end{enumerate}
\end{mainthm}

 Ramachandran and Rosenberg in \cite{ramachandran2025derived} show that  a genus 1 curve $C$ with Jacobian $J$ satisfies $\operatorname{D}^b(C)\simeq\operatorname{D}^b(J,\alpha_C)$ for a Brauer class $\alpha_C$ determined by the class of $C$ in $\cH^1(\bk,J)$. Hence any two genus 1 curves with isomorphic Jacobians are derived isogenous. From this, the implications (2) $\Rightarrow$ (3) $\Rightarrow $ (1) in the above statement follow easily. Our main input is the proof that derived isogenous genus 1 curves have the same Jacobian, which will be a consequence of the following:
 \begin{thm*}
 If $C$ and $C'$ are twisted derived equivalent genus 1 curves, then there is an isomorphism between their Jacobians $J(C)\simeq J(C')$.
  \end{thm*}
Our strategy is to leverage the main construction of \cite{olsson2025twisted}, where, starting from a twisted variety $(X,\alpha)$, the author provides a group algebraic space $\operatorname{O}(X,\alpha)$ which is a twisted derived invariant. More precisely, given a twisted equivalence $\Phi:\DXa\rightarrow \operatorname{D}^b(Y,\beta)$, Olsson constructs an isomorphism of group algebraic spaces 
$$R_\Phi:\operatorname{O}(X,\alpha)\rightarrow \operatorname{O}(Y,\beta).$$
When $X$ and $Y$ are genus 1 curves, using that the classes $\alpha$ and $\beta$ trivialize when passing to the algebraic closure of $\bk$, we construct, in Section \ref{Sec:OlssonTrivial}, isomorphisms
\[
\theta_\tau:J(X)\times \widehat{J(X)}\longrightarrow\operatorname{O}(X,\alpha),\quad \theta_{\tau'}:J(Y)\times \widehat{J(Y)}\longrightarrow\operatorname{O}(Y,\beta)
\]
such that the composition $\theta_{\tau'}^{-1}\circ R_\Phi\circ \theta_\tau$ is \emph{symplectic} in the sense of \cite{Polishchuk2003} (see also Subsection~\ref{subsec:symplectic} for more details). We can then conclude thanks to Proposition \ref{symplectic implies iso}, which shows that a symplectic isomorphism $E\times \hat E\rightarrow E'\times\hat{E'}$, where $E$ and $E'$ are elliptic curves, yields an isomorphism $E\simeq E'.$

\subsection*{Notation} We work over a perfect field $\bk$ with algebraic closure $\ol \bk$ and absolute Galois group $G_\bk=\operatorname{Gal}(\ol \bk/\bk)$. Given $X$ a reduced scheme over $\bk$ we denote by $\ol X$ its base extension over $\ol \bk$. For $\gamma\in G_\bk$ and an object or morphism $Z$ defined over $\ol \bk$ we write ${}^\gamma Z$ for its Galois conjugate.\par
A \emph{genus 1 curve} is a smooth projective geometrically connected curve of genus 1; its Jacobian $J(C)=\Pic^0_C$ is an elliptic curve.\par
Given an abelian variety $A$ over $\bk$, we denote by $\hat A$ its dual, and we identify $\hat{\hat A}=A$. Given a homomorphism $f\colon A\rightarrow B$ of abelian varieties over $\bk$, we denote by $\hat f\colon\hat B\rightarrow\hat A$ its dual. For an integer $m$ we denote by $[m]\colon A\rightarrow A$ the multiplication by $m$. For an elliptic curve $E$ we write $\lambda_E\colon E\to\hat E$ for its canonical principal polarization; recall that $\hat\lambda_E=\lambda_E$. \par
\subsection*{Acknowledgments}
The second and third authors would like to thanks D. Rydh for valuable mathematics conversation and M. Olsson for a very nice e-mail exchange.
The second author was partially supported by the grant number VR2023-03837 of the \emph{Vetenskaprådet} and by the grant number 2023.0395 of the Knut and Alice Wallenberg foundation. 
\subsection*{AI Usage} Claude Opus 5.5 was used to check for correand polish the exposition.

\section{Preliminary results and Background Material}

\subsection{Abelian varieties and Symplectic Isomorphisms}\label{subsec:symplectic}
Let $A$ and $B$ be abelian varieties over $\bk$. A morphism $f:A\times \hat A\rightarrow B\times \hat B$ can be represented by a matrix of morphisms:
$$
f=\begin{pmatrix}x&y\\z&w\end{pmatrix},
$$
with $x\colon A\to B$, $y\colon\hat A\to B$, $z\colon A\to\hat B$ and $w\colon\hat A\to\hat B$. We recall the following definition.
\begin{defi}\label{def:symplectic}
 An isomorphism $f\colon A\times\hat A\to B\times\hat B$ is said to be \emph{symplectic} (or \emph{isometric}) if
\begin{equation}\label{eq:inverseisometric}
    f^{-1}=\begin{pmatrix}\hw&-\hy\\-\hz&\hx\end{pmatrix}.
\end{equation}

\end{defi}
 The fact that two abelian varieties $A$ and $B$ are tied by a symplectic isomorphism $A\times \hat A\simeq B\times \hat B$ is very strong: over an algebraically closed field of characteristic not 2, this is equivalent to $A$ and $B$ being derived equivalent (\cite{orlov2002derived, polishchuk1996symplectic, Polishchuk2003}). Since for  elliptic curves derived equivalence implies isomorphism (\cite{antieau2017derived}), we directly get that elliptic curves tied by a symplectic isomorphism are isomorphic over (most) algebraically closed fields. The next result says that the same is true over any field.
 \begin{prop}\label{symplectic implies iso}
Let $\bk$ be a field and let $E_1,E_2$ be elliptic curves over $\bk$. If there is a symplectic isomorphism  {$f\colon E_1\times\hat E_1\to E_2\times\hat E_2$} over $\bk$, then $E_1$ and $E_2$ are isomorphic over $\bk$.
\end{prop}
\begin{proof}
Write $f$ in matrix form, that is $f=\begin{pmatrix}x&y\\z&w\end{pmatrix}$, and set $a:=\deg x$, $b:=\deg y$, $c:=\deg z$ and $d:=\deg w$.
Using the principal polarizations $\lambda_{E_i}\colon E_i\to\hat E_i$ we identify $\hat E_i$ with $E_i$. Then $x,y,z,w$ become elements of $\operatorname{Hom}_\bk(E_1,E_2)$, and the dual of a homomorphism $u\colon E_1\rightarrow E_2$ becomes its dual isogeny $\lambda_{E_1}^{-1}\circ \hat u\circ\lambda_{E_2}\colon E_2\rightarrow E_1$, which we still denote by $\hat u$; recall that $\hat u\circ u=[\deg u]$ and $u\circ\hat u=[\deg u]$.
 
Since $f$ is symplectic and the isomorphisms $\lambda_{E_i}$ are symmetric, $f^{-1}$ can be expressed as in \eqref{eq:inverseisometric}. Comparing the entries of the first row of $f\circ f^{-1}=\operatorname{Id}_{E_2\times\hat E_2}$ gives
\begin{equation}\label{eq:first row}
x\circ\hw-y\circ\hz=1_{E_2},\qquad y\circ\hx=x\circ\hy.
\end{equation}
If $x=0$, the first equation gives $-y\circ\hz=1_{E_2}$, so $y$ is an isomorphism and $E_1$ and $E_2$ are isomorphic. Similarly, if $z=0$, then $x\circ\hw=1_{E_2}$ and $x$ is an isomorphism. Therefore we can assume that $x,z\neq 0$, and thus $a,c\neq 0$.
 
As the dual of $y\circ \hx$ is $x\circ \hy$, the second equation in \eqref{eq:first row} implies that $y\circ\hx$ is self-dual, so it must be multiplication by an integer $n$. Indeed, the map $u\mapsto\hat u$ is the Rosati involution of the canonical principal polarization, whose fixed points in $\operatorname{End}(E_2)$ are the integers (see e.g.\ \cite[III.6.2 and \S III.9]{silverman2009arithmetic}).
 
By comparing degrees, we have $n^2=ab$; in particular $ab$ is a perfect square. Now we compose both sides of the second equation of \eqref{eq:first row} on the right with $x$. Using that $\hx\circ x=[a]$, we get
\begin{equation}\label{eq:ay=nx}
ay=nx.
\end{equation}
 
Now we turn our attention to the second row of $f\circ f^{-1}$, which gives the equations
\begin{equation}\label{second row}
z\circ\hw=w\circ\hz,\qquad w\circ\hx-z\circ\hy=1_{E_2}.
\end{equation}
As before, $z\circ \hw$ is self-dual, hence multiplication by an integer $m$ with $m^2=cd$. Composing on the right with $z$, and using $\hz\circ z=[c]$, we get
\begin{equation}\label{cw=mz}
cw=mz.
\end{equation}
 
By \eqref{eq:ay=nx} and \eqref{cw=mz} we have $y=\frac{n}{a}x$ in $\operatorname{Hom}(E_1,E_2)\otimes\mathbb{Q}$, and, taking duals, $\hw=\frac{m}{c}\hz$ in $\operatorname{Hom}(E_2,E_1)\otimes\mathbb{Q}$. Substituting into the first equation of \eqref{eq:first row}, we get the following equalities in $\operatorname{End}(E_2)\otimes\mathbb{Q}$:
\begin{equation}\label{end Q}
\begin{aligned}
1_{E_2}&=x\circ \left(\frac{m}{c}\hz\right)-\frac{n}{a}\,x\circ\hz\\
&=\left(\frac{m}{c}-\frac{n}{a}\right)x\circ\hz.
\end{aligned}
\end{equation}
This means that $x\circ\hz$ is a rational multiple $h$ of the identity. Comparing degrees, we get $h^2=ac$, and therefore $h$ is an integer. Now observe that, from \eqref{end Q}, we get

\[
1=\frac{hm}{c}-\frac{hn}{a}.
\]
Put $p:=hm/c$ and $q:=hn/a$. These are rational numbers with
\[
p^2=\frac{h^2m^2}{c^2}=\frac{ac\cdot cd}{c^2}=ad,
\qquad
q^2=\frac{h^2n^2}{a^2}=\frac{ac\cdot ab}{a^2}=bc,
\]
so $p$ and $q$ are integers. Since $p-q=1$, the integers $p$ and $q$ are relatively prime. Therefore the same is true for their squares $p^2=ad$ and $q^2=bc$. Any common divisor of $a$ and $b$ divides both $ad$ and $bc$. Thus we deduce that  $a$ and $b$ are coprime. As $ab=n^2$ is a perfect square, $a$ and $b$ are themselves perfect squares. Write $a=r^2$ with $r$ a positive integer.
 
Since $r^2=a$ divides $n^2=ab$, the integer $r$ divides $n$. Put $s:=n/r\in\Z$, so that $s^2=b$. Then \eqref{eq:ay=nx} reads $r^2y=rs\,x$. Since $\operatorname{Hom}_\bk(E_1,E_2)$ is torsion free and $r\neq0$, we get
\[
ry=sx.
\]
Moreover $r$ and $s$ are coprime, because $r^2=a$ and $s^2=b$ are. Choose $u,v\in\Z$ with $ur+vs=1$. Then
\[
x=(ur+vs)\,x=ur\,x+v\,sx=ur\,x+v\,ry=r\,(ux+vy).
\]
Again comparing degrees we have
$$a=\deg x=r^2\deg(ux+vy)=a\deg(ux+vy).$$ 
Therefore $\deg(ux+vy)=1$ because $a\neq0$. Hence $ux+vy\colon E_1\to E_2$ is an isomorphism defined over $\bk$, which concludes the proof.
\end{proof}
 We conclude by remarking that a variant of the proof of Proposition~\ref{symplectic implies iso} works for principally polarized abelian varieties $A$ and $B$ over $\bk$, with polarizations $\lambda_A$ and $\lambda_B$, whose only endomorphisms fixed by the Rosati involution are the integers, i.e.\ with $\operatorname{Hom}_\bk^{\mathrm{sym}}(A,\hat A)=\Z\lambda_A$ and $\operatorname{Hom}_\bk^{\mathrm{sym}}(B,\hat B)=\Z\lambda_B$. In that case the self-dual composites $y\circ\hx$ and $\hx\circ x$ are integers by hypothesis, the second one nonzero when $x\neq0$ by positivity of the Rosati involution. This shows that $x$ and $y$ are proportional, and one concludes by comparing degrees in the first equation of \eqref{eq:first row}, using $\dim A=\dim B$.
\subsection{Olsson groups}
In this paragraph we review the main construction in \cite{olsson2025twisted}. Before beginning, we observe that  the introduction of \cite{olsson2025twisted} is phrased over an algebraically closed field, but the constructions and proofs of the results we use are carried out over an arbitrary field; see \cite[3.1 and \S5]{olsson2025twisted}. \par
Let $X$ be a smooth projective geometrically connected variety over $\bk$ and let $\pi\colon\kx\to X$ be a $\IG_m$-gerbe with class $\alpha\in\operatorname{Br}(X)$. As in the introduction, we write $\DXa=\operatorname{D}^b(\kx)^{(1)}$ for the derived category of $1$-twisted sheaves. Let $\AutcX$  be the stack which to any $\bk$-scheme $T$ associates the groupoid of isomorphisms $\kx_T\to\kx_T$ inducing the identity on the stabilizer groups $\IG_m$.
By \cite{olsson2025twisted}, $\AutcX$ is a $\IG_m$-gerbe over a group algebraic space $\operatorname{Aut}_{\kx}$, locally of finite type over $\bk$. We set
\[
\operatorname{O}(X,\alpha):=\operatorname{Aut}^0_{\kx},
\]
the identity component, and we call it the \emph{Olsson group of $(X,\alpha)$}. Up to isomorphism, it only depends on $\alpha$, since gerbes associated to the same class are isomorphic.\par

By \cite[Theorem 1.1(2)]{olsson2025twisted} there is an exact sequence of group algebraic spaces
\begin{equation}\label{eq:olsson}
1\longrightarrow\Pic_{ X}\xrightarrow{\ \iota_{\mathscr X}\ }\operatorname{Aut}_{\mathscr X}\xrightarrow{\ r_{\mathscr X}\ }\Aut_{ X};
\end{equation}
more precisely, when the set is nonempty, the automorphisms of $\kx_T$ lying over a given automorphism of $X_T$ form a torsor under $\Pic_X(T)$ \cite[3.5 and Remark 3.6]{olsson2025twisted}. \par
 The map $\iota_{\kx}$ sends a line bundle $\kl$ to the automorphism $t_{\kl}$ of $\kx$ over $\operatorname{id}_X$ given by twisting by $\kl$. Explicitly, if $\xi$ is an object of $\kx$ over a $T$-point $x$ of $X$, then $t_{\kl}(\xi)=\xi\otimes x^*\kl$, where $\otimes$ denotes the action of $\IG_m$-torsors on $\kx$ given by its banding. (Depending on conventions, \cite[3.5]{olsson2025twisted} may use $\kl^{-1}$ instead of $\kl$; this plays no role below.)\par
The map $r_{\kx}$ sends $\sigma$ to the automorphism $\bar\sigma$ of $X$ that it induces. Since $X$ is the rigidification of $\kx$, every morphism from $\kx$ to an algebraic space factors uniquely through $\pi$. So $\bar\sigma$ is the unique automorphism of $X$ with $\pi\circ\sigma\simeq\bar\sigma\circ\pi$.\par

 Restricting to identity components we obtain an exact sequence
\begin{equation}\label{eq:olsson sequence}
0\rightarrow K_{X,\alpha}\rightarrow\operatorname{O}(X,\alpha)\rightarrow\Aut^0_X,\qquad K_{X,\alpha}:=\Pic_X\cap\operatorname{O}(X,\alpha).
\end{equation}

Since $\Pic^0_X$ is connected, $\Pic^0_X\subseteq K_{X,\alpha}$, but equality need not hold.  
On the other hand, $\operatorname{O}(X,\alpha)$ is of finite type, being the identity component of a group algebraic space locally of finite type over $\bk$. So $K_{X,\alpha}$ has finitely many connected components, and its image in $\operatorname{NS}(\ol X)=\Pic(\ol X)/\Pic^0(\ol X)$ is a finite group. Therefore
\[
\Pic^0_X\subseteq K_{X,\alpha}\subseteq\Pic^\tau_X .
\]
In particular $K_{X,\alpha}=\Pic^0_X$ whenever $\operatorname{NS}(\ol X)$ is torsion free, for instance when $X$ is a curve  or an abelian variety.\par
 \subsubsection{Olsson group for trivial Brauer classes}\label{subsec:trivialBC}

Let $\kx=B\IG_{m,X}$ be the trivial gerbe. Then
\[
\operatorname{D}^b(X)\xrightarrow{\ \sim\ }\operatorname{D}^b(B\IG_{m,X})^{(1)}=\DXa,\qquad \km\longmapsto\km\boxtimes\chi,
\]
where $\chi$ is the universal invertible sheaf on $B\IG_m$ (its pullback along the $T'$-point of $B\IG_m$ given by a line bundle $\mathcal N$ on $T'$ is $\mathcal N$). The stabilizers act on $\chi$ through the standard character. By \cite[Example 1.2]{olsson2025twisted} there is an isomorphism of stacks $\Aut_X\times\mathscr{P}ic_X\to\mathscr{A}ut_{B\IG_{m,X}}$ sending a pair $(a,\km)$, with $a\in\Aut_X(T)$ and $\km\in\mathscr{P}ic_X(T)$, to the automorphism $\sigma_{(a,\km)}$ which sends a $T$-point $(x,\mathcal N)$ of $X\times B\IG_m$ to $(a(x),\mathcal N\otimes x^*\km)$. %
Under this isomorphism the composition law reads
\begin{equation}\label{eq:composition}
\sigma_{(a,\km)}\circ\sigma_{(a',\km')}=\sigma_{(aa',\,\km'\otimes a'^*\km)},
\end{equation}
so that $\operatorname{Aut}_{B\IG_{m,X}}\simeq\Aut_X\ltimes\Pic_X$. If $\Pic^0_X$ is an abelian variety, then $\Aut^0_X$ acts trivially on $\Pic^0_X$ (a connected group scheme acts trivially on an abelian variety by group automorphisms), and therefore
\[
\operatorname{O}(X,0)\simeq\Aut^0_X\times\Pic^0_X
\]
as group algebraic spaces over $\bk$. In any case, the identity component of $\Aut_X\ltimes\Pic_X$ has underlying space $\Aut^0_X\times\Pic^0_X$. So, for the trivial gerbe the kernel in \eqref{eq:olsson sequence} is $K_{X,0}=\Pic^0_X$, and the map to $\Aut^0_X$ is surjective.
We will need the following two lemmas. The first one describes the action of $\operatorname{O}(X,0)(T)$ on $\operatorname{D}^b(X_T)$.
\begin{lemma}\label{action}
Let $(a,\km)\in\Aut_X(T)\times\mathscr{P}ic_X(T)$ and $\sigma=\sigma_{(a,\km)}$. Under the identification $\operatorname{D}^b(B\IG_{m,X_T})^{(1)}\simeq\operatorname{D}^b(X_T)$, the functor $\sigma^*$ is $\mathcal F\mapsto a^*\mathcal F\otimes\km$.
\end{lemma}
 
\begin{proof} For $a=\operatorname{id}$ this is simply \cite[Proposition 6.1]{olsson2025twisted} in the case of the trivial gerbe.\par
For the general case, let $p$ and $q$ be the projections from $X_T\times B\IG_m$ to first and second factor, respectively. By definition, $\sigma$ sends a $T$-point $(x,\mathcal N)$ to $(a(x),\mathcal N\otimes x^*\km)$. Hence $p\circ\sigma=a\circ p$. Moreover, the pullback of $(q\circ\sigma)^*\chi$ along $(x,\mathcal N)$ is $\mathcal N\otimes x^*\km$, which is the pullback of $q^*\chi\otimes p^*\km$ along the same point; since this holds functorially in $(x,\mathcal N)$, we get $(q\circ\sigma)^*\chi\simeq q^*\chi\otimes p^*\km$. Therefore
\[
\sigma^*(\mathcal F\boxtimes\chi)=(p\circ\sigma)^*\mathcal F\otimes(q\circ\sigma)^*\chi\simeq p^*a^*\mathcal F\otimes q^*\chi\otimes p^*\km=(a^*\mathcal F\otimes\km)\boxtimes\chi.
\]

\end{proof}
 The second lemma is about the kernel of the action on $\bk$-valued points when $\bk$ is an algebraically closed field.
\begin{lemma}\label{faithful}
Assume that $\bk$ is algebraically closed. If $(a,\km)\in\Aut_X(\bk)\times\Pic_X(\bk)$ and $a^*(-)\otimes\km\simeq\operatorname{id}$ as functors on $\operatorname{D}^b(X)$, then $a=\operatorname{id}_X$ and $\km\simeq\mathcal O_X$.
\end{lemma}
 
\begin{proof}
This is a special case of \cite[Remark 1.8]{olsson2025twisted}, which deduces the faithfulness of the action of $\AutcX(\bk)$ on $\DXa$ from the uniqueness of Fourier--Mukai kernels \cite[Theorem 1.1]{canonaco2007twisted}. Directly: evaluating on the skyscraper sheaves $\bk(x)$ gives $a^{-1}(x)=x$ for all closed points $x$, so $a=\operatorname{id}_X$, and evaluating on $\mathcal O_X$ gives $\km\simeq\mathcal O_X$.
\end{proof}
\subsubsection{Twisted Rouquier isomorphisms}\label{subsec:rouquier}
The main result of \cite{olsson2025twisted} is that the Olsson group is a twisted derived invariant. Let $Y$ be another smooth projective geometrically connected variety over $\bk$, let $\mathscr Y\to Y$ be a $\IG_m$-gerbe with class $\beta\in\operatorname{Br}(Y)$, and let $\Phi\colon\DXa\to\operatorname{D}^b(Y,\beta)$ be an exact equivalence. By \cite[Theorem 1.1]{canonaco2007twisted}, $\Phi$ is of Fourier--Mukai type, and by \cite[Theorems 1.7 and 1.9]{olsson2025twisted} it induces an isomorphism of group algebraic spaces
\[
R_\Phi\colon\operatorname{O}(X,\alpha)\longrightarrow\operatorname{O}(Y,\beta),
\]
which we call the \emph{twisted Rouquier isomorphism} associated to $\Phi$.

Let an automorphism $\sigma$ of a gerbe act on the twisted derived category by $\sigma^*$, as in \cite[Example 1.5]{olsson2025twisted}. Then, on $\ol\bk$-points, $R_\Phi$ is given by conjugation \cite[Remark 7.2]{olsson2025twisted}:
\begin{equation}\label{eq:conjugation}
\Phi_{\ol\bk}\circ\sigma^*\circ\Phi_{\ol\bk}^{-1}\simeq R_\Phi(\sigma)^*.
\end{equation}
When $\alpha$ and $\beta$ are trivial and $\operatorname{O}(X,0)$ and $\operatorname{O}(Y,0)$ are reduced, then the twisted Roquier isomorphism is the isomorphism constructed by Roquier in  \cite{rouquier2011automorphismes}.
In particular, when $X$ and $Y$ are abelian varieties with an equivalence $\Phi\colon\operatorname{D}^b(X)\to\operatorname{D}^b(Y)$, we have $\operatorname{O}(A,0)=X\times\hat X$ and $\operatorname{O}(B,0)=Y\times\hat Y$, up to some negligible conventions discussed in Remark \ref{rmk:conventions}. Then the twisted Roquier isomorphism is exactly  the isomorphism $f_\Phi$ constructed by Orlov \cite{orlov2002derived}. In particular, $R_\Phi$ is symplectic by \cite[Corollary 2.13]{orlov2002derived}.
\begin{rmk}[Conventions]\label{rmk:conventions}
The comparison above depends on two conventions: whether an automorphism $a$ acts on $\operatorname{D}^b(X)$ through $a^*$ or through $a_*=(a^{-1})^*$, and whether a line bundle $\kl$ acts through $-\otimes\kl$ or through $-\otimes\kl^{-1}$ (see the parenthetical remark after \eqref{eq:olsson}). These choices do not affect symplecticity. Indeed, let $A$ and $B$ abelian varieties with $\Phi:\operatorname{D}^b(A)\simeq \operatorname{D}^b(B)$ an equivalence between their untwisted derived categories. The isomorphism $f_\Phi:A\times \hat{A}\rightarrow B\times \hat{B}$ constructed by Orlov, is such that, over $\ol{\bk}$, we have $f_\Phi(a,\xi)=(b,\eta)$ if, and only if
\begin{equation}\label{eq:orlov}
\Phi\circ\big(t_a^*(-)\otimes\mathcal P_\xi\big)\circ\Phi^{-1}\simeq t_b^*(-)\otimes\mathcal P_\eta,
\end{equation} 
where $t_a$ denotes the translation by $a$ and $\mathcal P_\xi$ denotes the line bundle in $\Pic^0(A)$ corresponding to $\xi\in\hat A(\ol\bk)$. Now, let $f'$ be the isomorphism characterized by the analogue of \eqref{eq:orlov} in which $t_a^*$ and $t_b^*$ are replaced by $t_{a*}=t_{-a}^*$ and $t_{b*}=t_{-b}^*$. Then $f'(a,\xi)=(b,\eta)$ if and only if $f(-a,\xi)=(-b,\eta)$, that is,
\[
f'=\epsilon_B\circ f\circ\epsilon_A,\qquad\text{where }\epsilon_A:=\begin{pmatrix}-1&0\\0&1\end{pmatrix}\colon(a,\xi)\longmapsto(-a,\xi).
\]
Similarly, replacing $\mathcal P_\xi$ and $\mathcal P_\eta$ by $\mathcal P_\xi^{-1}=\mathcal P_{-\xi}$ and $\mathcal P_\eta^{-1}=\mathcal P_{-\eta}$ replaces $f$ by $\epsilon'_B\circ f\circ\epsilon'_A$, where $\epsilon'_A:=-\epsilon_A$; which  is again $\epsilon_B\circ f\circ\epsilon_A$ (and changing both conventions leaves $f$ unchanged). In matrix form $f"$ is represented as follows
\[
\epsilon_B\circ\begin{pmatrix}x&y\\z&w\end{pmatrix}\circ\epsilon_A=\begin{pmatrix}x&-y\\-z&w\end{pmatrix}.
\]
If $f$ is symplectic, then by \eqref{eq:inverseisometric} the inverse of this isomorphism is
\[
\epsilon_A\circ f^{-1}\circ\epsilon_B=\epsilon_A\circ\begin{pmatrix}\hw&-\hy\\-\hz&\hx\end{pmatrix}\circ\epsilon_B=\begin{pmatrix}\hw&\hy\\\hz&\hx\end{pmatrix},
\]
which is exactly the matrix that Definition~\ref{def:symplectic} prescribes for $\begin{pmatrix}x&-y\\-z&w\end{pmatrix}$, because $\widehat{-u}=-\hat u$. Hence $\epsilon_B\circ f\circ\epsilon_A$ is symplectic whenever $f$ is, and we may ignore these conventions. \qed 
\end{rmk}

We now describe two twisted Rouquier isomorphisms that play a central role in Section~3.

\begin{ex}[Conjugation by an isomorphism of gerbes]\label{ex:Rtau}
Let $\kx$ and $\kx'$ be $\IG_m$-gerbes over $X$, and let $\tau\colon\kx\to\kx'$ be an isomorphism of $\IG_m$-gerbes over $X$, that is, an isomorphism of stacks which commutes with the projections to $X$ and induces the identity on stabilizer groups. Then $\tau^*\colon\operatorname{D}^b(\kx')^{(1)}\to\operatorname{D}^b(\kx)^{(1)}$ is an equivalence and thus induces an isomorphism $R_{\tau^*}\colon\operatorname{O}(X,\beta)\rightarrow \operatorname{O}(X,\alpha)$, where $\alpha$ and $\beta$ are the Brauer classes associated to $\mathscr{X}\rightarrow X$ and $\mathscr{X'}\rightarrow X$ (of course $\alpha=\beta$, since the two gerbes are isomorphic). By \eqref{eq:conjugation} and the faithfulness of the action \cite[Remark 1.8]{olsson2025twisted}, we see that $R_{\tau^*}$ is the isomorphism induced after rigidification by the isomorphism of group stacks $\mathscr{A}ut_{\kx'}\to\mathscr{A}ut_{\kx}$ given by conjugation by $\tau$: $\sigma\mapsto\tau^{-1}\circ\sigma\circ\tau$.\par
To lighten the notation we will, in the sequel, denote $R_{\tau^*}$ simply by $R_\tau$.\par
We see that  $R_\tau$ behaves well with respect to the exact sequences \eqref{eq:olsson} for $\kx'$ and $\kx$, i.e.\ the diagram
\begin{equation}\label{eq:ctau}
\begin{tikzcd}[column sep=large]
1\arrow[r] & \Pic_{X}\arrow[r,"\iota_{\kx'}"]\arrow[d,equal] & \operatorname{Aut}_{\kx'}\arrow[r,"r_{\kx'}"]\arrow[d,"R_\tau"] & \Aut_{X}\arrow[d,equal]\\
1\arrow[r] & \Pic_{X}\arrow[r,"\iota_{\kx}"] & \operatorname{Aut}_{\kx}\arrow[r,"r_{\kx}"] & \Aut_{X}
\end{tikzcd}
\end{equation}
commutes. Let $\pi\colon\kx\to X$ and $\pi'\colon\kx'\to X$ be the structure maps. Saying that $\tau$ is an isomorphism of gerbes over $X$ means that $\pi'\circ\tau\simeq\pi$. Hence, for every automorphism $\sigma$ of $\kx'$,
\[
\pi\circ(\tau^{-1}\sigma\tau)\simeq\pi'\circ\sigma\circ\tau\simeq\bar\sigma\circ\pi'\circ\tau\simeq\bar\sigma\circ\pi.
\]
By the uniqueness in the definition of $r$, this gives $r_{\kx}(\tau^{-1}\circ\sigma\circ \tau)=\bar\sigma=r_{\kx'}(\sigma)$, so the right square commutes. For the left square, observe that $\tau$ induces the identity on stabilizers, so it is compatible with the action of $\IG_m$-torsors: $\tau(\xi\otimes P)\simeq\tau(\xi)\otimes P$, functorially in the object $\xi$ and in the $\IG_m$-torsor $P$. Taking $P=x^*\kl$ gives $\tau\circ t_{\kl}\simeq t_{\kl}\circ\tau$, that is, $R_\tau(\iota_{\kx'}(\kl))=\iota_{\kx}(\kl)$.\par

The case of main interest for us will be when, over $\ol\bk$, the gerbe $\kx_{\ol\bk}\rightarrow \ol{X}$ is associated to the trivial Brauer class and $\tau\colon\kx_{\ol\bk}\rightarrow B\IG_{m,\ol X}$ is a trivialization.  Applying the above over $\ol\bk$ and passing to identity components, by Subsection~\ref{subsec:trivialBC} the top row of \eqref{eq:ctau} restricts to the split exact sequence $0\to\Pic^0_{\ol X}\to\operatorname{O}(X,0)_{\ol \bk}\to\Aut^0_{\ol X}\to0$. We therefore obtain an isomorphism of group algebraic spaces $R_\tau\colon\operatorname{O}(X,0)_{\ol\bk}\to\operatorname{O}(X,\alpha)_{\ol\bk}$ fitting in the commutative diagram
\begin{equation}\label{eq:ctau0}
\begin{tikzcd}[column sep=small]
0\arrow[r] & \Pic^0_{\ol X}\arrow[r]\arrow[d,equal] & \operatorname{O}(X,0)_{\ol \bk}\arrow[r]\arrow[d,"R_\tau"] & \Aut^0_{\ol X}\arrow[r]\arrow[d,equal] & 0\\
0\arrow[r] & \Pic^0_{\ol X}\arrow[r] & \operatorname{O}(X,\alpha)_{\ol \bk}\arrow[r] & \Aut^0_{\ol X}\arrow[r] & 0
\end{tikzcd}
\end{equation}
The bottom row is exact, because the top row is and $R_\tau$ is an isomorphism. When $\Pic^0_X$ is an abelian variety, $\operatorname{O}(X,0)=\Aut^0_X\times\Pic^0_X$.
\end{ex}

\begin{ex}[The automorphisms $\operatorname{Ad}(\mathcal Q)$]\label{ex:Ad}
Let $\mathcal Q$ be a line bundle on $X$. Then $\sigma_{(\operatorname{id},\mathcal Q)}$ is an automorphism of the trivial gerbe $B\IG_{m,X}$, with inverse $\sigma_{(\operatorname{id},\mathcal Q)}^{-1}=\sigma_{(\operatorname{id},\mathcal Q^{-1})}$ by \eqref{eq:composition}. Applying Example~\ref{ex:Rtau} to $\tau=\sigma_{(\operatorname{id},\mathcal Q^{-1})}$, we obtain the automorphism
\[
\operatorname{Ad}(\mathcal Q):=R_{\sigma_{(\operatorname{id},\mathcal Q^{-1})}}\colon\rho\longmapsto\sigma_{(\operatorname{id},\mathcal Q)}\circ\rho\circ\sigma_{(\operatorname{id},\mathcal Q)}^{-1}
\]
of $\operatorname{Aut}_{B\IG_{m,X}}$, that is, conjugation by $\sigma_{(\operatorname{id},\mathcal Q)}$. By Example~\ref{ex:Rtau} and Lemma~\ref{action}, it is the twisted Rouquier isomorphism of the autoequivalence $\sigma_{(\operatorname{id},\mathcal Q^{-1})}^*=-\otimes\mathcal Q^{-1}$ of $\operatorname{D}^b(X)$. Note that $\operatorname{Ad}(\mathcal Q)$ is not an automorphism of $B\IG_{m,X}$ itself: it acts on automorphisms of $B\IG_{m,X}$. It preserves the identity component $\operatorname{O}(X,0)$, and we denote the induced automorphism again by $\operatorname{Ad}(\mathcal Q)$.

By \eqref{eq:composition} we have $\operatorname{Ad}(\mathcal Q\otimes\mathcal Q')=\operatorname{Ad}(\mathcal Q)\circ\operatorname{Ad}(\mathcal Q')$, and
\begin{equation}\label{eq:Ad}
\operatorname{Ad}(\mathcal Q)(\sigma_{(a,\kl)})=\sigma_{(a,\ \kl\otimes a^*\mathcal Q\otimes\mathcal Q^{-1})}.
\end{equation}
Assume now that $\Pic^0_X$ is an abelian variety, so that $\operatorname{O}(X,0)=\Aut^0_X\times\Pic^0_X$. Then $\operatorname{Ad}(\mathcal Q)$ acts on it by $(a,\kl)\mapsto(a,\kl\otimes a^*\mathcal Q\otimes\mathcal Q^{-1})$. If $\mathcal Q\in\Pic^0(X)$, this automorphism is the identity: the morphism $a\mapsto a^*\mathcal Q\otimes\mathcal Q^{-1}$ from $\Aut^0_{X}$ to $\Pic^0_{X}$ is trivial, because $\Aut^0_X$ acts trivially on the abelian variety $\Pic^0_X$. Consequently, on $\Aut^0_X\times\Pic^0_X$ the automorphism $\operatorname{Ad}(\mathcal Q)$ only depends on the class of $\mathcal Q$ in $\operatorname{NS}(X)$. In Section~3 we use this over $\ol\bk$.
\end{ex}

\section{Olsson groups for Algebraic Brauer classes}\label{Sec:OlssonTrivial}

The main result of this section is the following
\begin{prop}\label{split2}
    If $C$ is a genus 1 curve and $\alpha\in \operatorname{Br}(C)$, then 
    $$\operatorname{O}(C,\alpha)\simeq J(C)\times \hat J(C).$$
\end{prop}

We say that $\alpha\in\operatorname{Br}(X)$ is \emph{algebraic} when it sits in the 
kernel of the map $\operatorname{Br}(X)\to \operatorname{Br}(\ol X)$. Note that every class $\alpha\in\operatorname{Br}(C)$ becomes trivial over $\ol\bk$: since $\ol C$ is a smooth projective curve over an algebraically closed field, $\operatorname{Br}(\ol C)=0$ by Tsen's theorem. Thus every $\alpha\in\operatorname{Br}(C)$ is algebraic.\par

We will put ourselves in a much more general situation and work with $X$, a smooth geometrically connected projective variety over a perfect field $\bk$ such that $\Pic^0_X$ is an abelian variety. The idea is that, when $\alpha$ is algebraic, over $\ol \bk$ there is an isomorphism
\begin{equation*}
    \operatorname{O}(X,\alpha)_{\ol{\bk}}\simeq \left(\Pic^0_X\times\operatorname{Aut}^0_X\right)_{\ol{\bk}},
\end{equation*}
since $\ol  \alpha$ is trivial by hypothesis and $\operatorname{Pic}^0_X$ is an abelian variety (see Subsection \ref{subsec:trivialBC}). So we want to define a morphism $\operatorname{Aut}^0_X\times\operatorname{Pic}^0_X\to \operatorname{O}(X,\alpha)$ over $\bk$ whose base change to $\ol \bk$ is the isomorphism $R_\tau$ of Example~\ref{ex:Rtau}. \par

If $\alpha$ is algebraic, then $\kx_{\ol \bk}$ is a trivial gerbe, so there is an isomorphism of $\IG_m$-gerbes
\[
\tau\colon\kx_{\ol \bk}\xrightarrow{\ \sim\ }B\IG_{m,\ol X},
\]
which, by a limit argument, is defined over some finite Galois extension $\bk'/\bk$. \par
Let $\gamma\in G_\bk$. The projections $\kx_{\ol \bk}\to\ol X$ and $B\IG_{m,\ol X}\to\ol X$ are base changes of morphisms defined over $\bk$, hence Galois equivariant. So the conjugate ${}^\gamma\tau$ again commutes with the projections to $\ol X$ and induces the identity on stabilizers. Consider the composite
\[
\rho_\gamma:={}^\gamma\tau\circ\tau^{-1}\colon B\IG_{m,\ol X}\to B\IG_{m,\ol X}.
\]
Clearly, it induces the identity on stabilizers, so it is an object of $\mathscr{A}ut_{B\IG_{m,\ol X}}(\ol \bk)$. In addition, it commutes with the projection to $\ol X$, so its image under the map $\mathscr{A}ut_{B\IG_{m,\ol X}}\to\Aut_{\ol X}$ of \cite[Theorem 1.1(2)]{olsson2025twisted} is $\operatorname{id}_{\ol X}$. Hence, by \cite[Example 1.2]{olsson2025twisted} (see Subsection~\ref{subsec:trivialBC}), it is of the form $\sigma_{(\operatorname{id},\km_\gamma)}$ for a line bundle $\km_\gamma$ on $\ol X$. By Lemma \ref{faithful} this line bundle is  unique up to isomorphism.

\begin{lemma}\label{lem:cocycle}
In the notation above we have the following:
\begin{enumerate}
\item for $\gamma\in\operatorname{Gal}(\ol \bk/\bk')$, $\km_\gamma\simeq\mathcal O_{\ol X}$;
\item for all $\gamma,\delta\in G_\bk$, $\km_{\gamma\delta}\simeq\km_\gamma\otimes{}^\gamma\km_\delta$;

\item every other trivialization of $\kx_{\ol \bk}$ is of the form $\tau'=\sigma_{(\operatorname{id},\mathcal Q)}\circ\tau$ with $\mathcal Q\in\Pic(\ol X)$. In addition, if we denote by $\km'_\gamma$ the line bundle associated to ${}^{\gamma}\tau'\circ\tau'^{-1}$ then $$\km'_\gamma\simeq\km_\gamma\otimes{}^\gamma\mathcal Q\otimes\mathcal Q^{-1}.$$
\end{enumerate}
\end{lemma}
 
\begin{proof} The first claim is clear since $\tau$ is defined over $\bk'$, and thus left fixed by conjugation by $\gamma$. We turn our attention to (2). Since the isomorphism $\Aut_X\times\mathscr{P}ic_X\simeq\mathscr{A}ut_{B\IG_{m,X}}$ is defined over $\bk$, we have that 
\begin{align*}
    {}^{\gamma\delta}\tau\circ\tau^{-1}&={}^{\gamma}({}^\delta\tau)\circ\tau^{-1}\\
    &={}^{\gamma}({}^\delta\tau\circ \tau^{-1})\circ{}^\gamma\tau\circ\tau^{-1}\\
    &=\sigma_{(\operatorname{id},{}^\gamma\km_\delta)}\circ \sigma_{(\operatorname{id},\km_\gamma)}.
\end{align*}
The result follows from \eqref{eq:composition}.\par
 
For (3), observe that the composite $\tau'\circ\tau^{-1}$ is an automorphism of $B\IG_{m,\ol X}$ over the identity inducing the identity on stabilizers. Thus 
\[\tau'\circ\tau^{-1}=\sigma_{(\operatorname{id},\mathcal{Q})}\]
for some line bundle $\mathcal{Q}$. Hence $\tau'$ has the stated form. Moreover $${}^\gamma\tau'\circ\tau'^{-1}=\sigma_{(\operatorname{id},{}^\gamma\mathcal Q)}\circ\sigma_{(\operatorname{id},\km_\gamma)}\circ\sigma_{(\operatorname{id},\mathcal Q^{-1})}=\sigma_{(\operatorname{id},\km_\gamma\otimes{}^\gamma\mathcal Q\otimes\mathcal Q^{-1})},$$ 
where the last step is \eqref{eq:composition}, and the statement is proven.
 
\end{proof}
By Lemma~\ref{lem:cocycle}(2), the map $\gamma\mapsto [\km_\gamma]\in \operatorname{NS}(\overline{X})$ is a $1$-cocycle, so it defines a class in $\cH^1(\bk,\operatorname{NS}(\overline{X}))$. By Lemma~\ref{lem:cocycle}(3), replacing $\tau$ by another trivialization $\sigma_{(\operatorname{id},\mathcal Q)}\circ\tau$ changes this cocycle by $\gamma\mapsto[{}^\gamma\mathcal Q\otimes\mathcal Q^{-1}]={}^\gamma[\mathcal Q]-[\mathcal Q]$, which is a coboundary; hence the class does not depend on the choice of $\tau$. It does not depend on the gerbe $\kx$ either, but only on $\alpha$. In fact, given another gerbe $\kx'$ representing $\alpha$ and an isomorphism $i\colon\kx'\rightarrow\kx$ over $X$, we have 
$${}^\gamma(\tau\circ i)\circ (\tau\circ i)^{-1}={
}^\gamma\tau\circ \tau^{-1}.
 $$
 Thus we can define the following
 \begin{defi}
     For $\alpha \in\operatorname{Br}(X)$ algebraic we define the class
     $$
     \operatorname{ob}(\alpha):=\big[\gamma\mapsto\km_\gamma \bmod \Pic^0(\ol X)\big]\in\cH^1(\bk,\operatorname{NS}(\ol X)).
     $$
 \end{defi}

\begin{rmk}\label{rem:ob}

If $\operatorname{NS}(\ol X)$ is torsion free with trivial Galois action, then
\[
\cH^1(\bk,\operatorname{NS}(\ol X))=
0,
\]
hence $\operatorname{ob}(\alpha)=0$. This applies to genus 1 curves, for which $\operatorname{NS}(\ol C)=\Z$ via the degree.

\end{rmk}
From now on, until the end of this section, we work under the following assumptions.
\begin{setting}\label{assumption}
    Let $X$ be a smooth projective geometrically connected variety over $\bk$ such that $\Pic^0_X$ is an abelian variety, and let $\alpha\in\operatorname{Br}(X)$ be algebraic.
\end{setting}

Let $\tau\colon\kx_{\ol \bk}\to B\IG_{m,\ol X}$ be a trivialization. By Example~\ref{ex:Rtau}, conjugation by $\tau$ gives an isomorphism of group algebraic spaces over $\ol\bk$
\[
R_\tau\colon(\Aut^0_X\times\Pic^0_X)_{\ol \bk}\longrightarrow\operatorname{O}(X,\alpha)_{\ol \bk},
\]
fitting in the commutative diagram \eqref{eq:ctau0} with exact rows. Here we use that $\operatorname{O}(X,0)\simeq\Aut^0_X\times\Pic^0_X$, because $\Pic^0_X$ is an abelian variety (Subsection~\ref{subsec:trivialBC}). For a line bundle $\mathcal Q$ on $\ol X$ we also use the automorphism $\operatorname{Ad}(\mathcal Q)$ of $(\Aut^0_X\times\Pic^0_X)_{\ol \bk}$ of Example~\ref{ex:Ad}. It is the identity if $\mathcal Q\in\Pic^0(\ol X)$, and in general it only depends on the class of $\mathcal Q$ in $\operatorname{NS}(\ol X)$.
\begin{lemma}\label{lem:trivializations}
Under the assumption in Setting \ref{assumption}, we have the following:
\begin{enumerate}
\item If $\operatorname{ob}(\alpha)=0$, there is a trivialization $\tau\colon\kx_{\ol \bk}\to B\IG_{m,\ol X}$ with $\km_\gamma\in\Pic^0(\ol X)$ for all $\gamma\in G_\bk$.
\item Let $\tau$ and $\tau'$ be two trivializations of $\kx_{\ol\bk}$, and write $\tau'=\sigma_{(\operatorname{id},\mathcal Q)}\circ\tau$ with $\mathcal{Q}\in\Pic(\ol X)$ as in Lemma~\ref{lem:cocycle}(3). Then $R_{\tau'}=R_\tau\circ\operatorname{Ad}(\mathcal Q)^{-1}$. If moreover $\operatorname{ob}(\alpha)=0$ and  both $\tau$ and $\tau'$ are as in (1), then the class of $\mathcal Q$ in $\operatorname{NS}(\ol X)$ is Galois invariant.
\end{enumerate}
\end{lemma}

\begin{proof}
(1) Let $\tau_0$ be any trivialization, with associated line bundles $\km^0_\gamma$. Since $\operatorname{ob}(\alpha)=0$, there is $\nu\in\operatorname{NS}(\ol X)$ such that the class of $\km^0_\gamma$ is ${}^\gamma\nu-\nu$ for all $\gamma\in G_\bk$. Let $\mathcal Q\in\Pic(\ol X)$ be a lift of $-\nu$ and put $\tau=\sigma_{(\operatorname{id},\mathcal Q)}\circ\tau_0$. Denote the associated line bundles by $\km_{\gamma}$. By Lemma~\ref{lem:cocycle}(3), the class of $\km_\gamma$ in $\operatorname{NS}(\ol X)$ is $({}^\gamma\nu-\nu)-({}^\gamma\nu-\nu)=0$, that is, $\km_\gamma\in\Pic^0(\ol X)$.

(2) For every $\rho$ we have
\[
R_{\tau'}(\rho)=\tau^{-1}\circ\sigma_{(\operatorname{id},\mathcal Q)}^{-1}\circ\rho\circ\sigma_{(\operatorname{id},\mathcal Q)}\circ\tau=R_\tau\big(\operatorname{Ad}(\mathcal Q)^{-1}(\rho)\big).
\]
If both $\tau$ and $\tau'$ are as in (1), Lemma~\ref{lem:cocycle}(3) gives ${}^\gamma\mathcal Q\otimes\mathcal Q^{-1}\simeq\km'_\gamma\otimes\km_\gamma^{-1}\in\Pic^0(\ol X)$ for all $\gamma\in G_\bk$, so the class of $\mathcal Q$ in $\operatorname{NS}(\ol X)$ is Galois invariant.
\end{proof}

The next result shows that $\operatorname{O}(X,\alpha)\simeq\operatorname{Aut}^0_X\times\operatorname{Pic}^0_X$ when the obstruction $\operatorname{ob}(\alpha)$ vanishes. Proposition~\ref{split2} then follows immediately from Remark~\ref{rem:ob} and the classical facts that, for a genus 1 curve $C$, we have $\operatorname{Pic}^0_C\simeq J(C)$ and $\operatorname{Aut}_C^0\simeq \widehat{J(C)}$.
\begin{prop}\label{split}
Under the hypothesis of Setting \ref{assumption}, assume in addition that $\operatorname{ob}(\alpha)=0$, and let $\tau\colon\kx_{\ol \bk}\to B\IG_{m,\ol X}$ be a trivialization as in Lemma~\ref{lem:trivializations}(1). Then $R_\tau$ descends to an isomorphism of group algebraic spaces over $\bk$
\[
\theta_\tau\colon\Aut^0_X\times\Pic^0_X\xrightarrow{\ \sim\ }\operatorname{O}(X,\alpha),
\]
compatible with the inclusions of $\Pic^0_X$ and with the maps to $\Aut^0_X$; that is, $\theta_\tau$ fits in the analog of \eqref{eq:ctau0} over $\bk$. In particular, in \eqref{eq:olsson sequence} the map $\operatorname{O}(X,\alpha)\to\Aut^0_X$ is surjective and $K_{X,\alpha}=\Pic^0_X$, so that $\operatorname{O}(X,\alpha)$ is a split extension of $\Aut^0_X$ by $\Pic^0_X$.
\end{prop}

\begin{proof}
Once we know that $R_\tau$ descends to an isomorphism $\theta_\tau$ of group algebraic spaces over $\bk$, we are done. In fact, the commutativity of the analog of \eqref{eq:ctau0} over $\bk$ can be checked over $\ol \bk$, where it is \eqref{eq:ctau0} itself. The surjectivity of $\operatorname{O}(X,\alpha)\to\Aut^0_X$ and the equality $K_{X,\alpha}=\Pic^0_X$ can also be checked over $\ol \bk$, where they follow from the exactness of the bottom row of \eqref{eq:ctau0}. Finally, $\theta_\tau$ itself gives the splitting.\par
To construct $\theta_\tau$, we use Galois descent. We first observe that the trivialization $\tau$ is defined over some finite Galois extension $\bk'$ of $\bk$, and thus so is $R_\tau$.\par
Now, for $\gamma\in G_\bk$, we have ${}^\gamma(R_\tau)=R_{{}^\gamma\tau}$. In addition, by definition of $\km_\gamma$ we have $${}^\gamma\tau=\rho_\gamma\circ\tau=\sigma_{(\operatorname{id},\km_\gamma)}\circ\tau.$$
Hence ${}^\gamma(R_\tau)=R_\tau\circ\operatorname{Ad}(\km_\gamma)^{-1}$ by Lemma~\ref{lem:trivializations}(2).\par
Since $\tau$ is as in Lemma~\ref{lem:trivializations}(1), we have $\km_\gamma\in\Pic^0(\ol X)$ for all $\gamma\in G_\bk$. Hence, as observed after \eqref{eq:Ad}, the automorphism $\operatorname{Ad}(\km_\gamma)$ of $(\Aut^0_X\times\Pic^0_X)_{\ol \bk}$ is the identity. Therefore $R_\tau$ is invariant under the action of $G_\bk$, which factors through $\operatorname{Gal}(\bk'/\bk)$. So $R_\tau$ descends to a morphism $\theta_\tau$ over $\bk$. The same applies to its inverse, and thus $\theta_\tau$ is an isomorphism. It is a homomorphism of group algebraic spaces, since this can be checked over $\ol\bk$, where it holds for $R_\tau$.
\end{proof}
We now specialize to the case relevant for the main result: torsors under abelian varieties. To this aim, we let $X$ be a torsor under an abelian variety $A$ over $\bk$. Our main example will be a genus 1 curve, which is a torsor under its Jacobian. The action of $A$ on $X$ by translations identifies $\operatorname{Aut}^0_X$ with $A$ \cite[Lemma A.3]{honigs2018derived}. Moreover, $\Pic^0_X$ is canonically isomorphic to $\Pic^0_A=\hat A$ \cite[XIII, 1.1]{raynaud1970faisceaux} (see also \cite[\S4]{ramachandran2025derived}). \par
In particular $\Pic^0_X$ is an abelian variety. If $\alpha\in\operatorname{Br}(X)$ is an algebraic class, we are under the assumptions of Setting~\ref{assumption}. Thus, if $\operatorname{ob}(\alpha)=0$, Proposition~\ref{split} gives isomorphisms
\[
\theta_\tau\colon A\times\hat A\xrightarrow{\ \sim\ }\operatorname{O}(X,\alpha).
\]

\section{Proof of the main theorem}\label{sec:symplectic}
In this section we prove the Main Theorem from the introduction.\par
To this aim, let $C$ and $C'$ be two genus 1 curves over $\bk$ with Jacobians $J$ and $J'$ respectively. If  (2) hold, then they are torsors over the same abelian variety $J$, and there are two classes $\alpha$ and $\alpha'$ in $\cH^1(\bk,J)$ such that $[C]=\alpha$ and $[C']=\alpha'$. Regarding them as Brauer classes on $J$ via $\operatorname{Br}(J)\simeq\operatorname{Br}(\bk)\oplus\cH^1(\bk,J)$, by \cite{ramachandran2025derived} we have twisted equivalences
$$
\operatorname{D}^b(C)\simeq \operatorname{D}^b(J,\alpha),\quad \operatorname{D}^b(J,\alpha')\simeq \operatorname{D}^b(C').
$$
We deduce that (3) holds, which  implies directly (1). To conclude the proof we need to show that two derived isogenous genus 1 curves have isomorphic Jacobians.  Suppose that $C$ and $C'$ are derived isogenous. Then there is a sequence of varieties $C=X_0,\ldots X_n=C'$  and with twisted equivalences
$$
\operatorname{D}^b(X_i,\beta_i)\simeq\operatorname{D}^b(X_{i+1},\alpha_{i+1}).
$$
Since twisted derived equivalence preserve dimension and Kodaira dimension (or order of the canonical bundle), we can show by induction that $X_i$ are curves of genus 1. Therefore, the last implication of the main theorem will be a consequence of Proposition~\ref{symplectic implies iso} and the following result.

\begin{prop}\label{Twisted implies symplectic}
Let $C$ and $C'$ be as before, and let $\alpha\in\operatorname{Br}(C)$ and $\beta\in\operatorname{Br}(C')$. If there is an exact equivalence $\operatorname{D}^b(C,\alpha)\simeq\operatorname{D}^b(C',\beta)$, then there is a symplectic isomorphism $J\times\hat J\simeq J'\times\hat{J'}$ defined over $\bk$.
\end{prop}

\begin{proof}
Let $\Phi\colon\operatorname{D}^b(C,\alpha)\to\operatorname{D}^b(C',\beta)$ be an equivalence, and let $R_\Phi\colon\operatorname{O}(C,\alpha)\to\operatorname{O}(C',\beta)$ be its twisted Rouquier isomorphism (Subsection~\ref{subsec:rouquier}). Let $\kx\to C$ and $\mathscr Y\to C'$ be the gerbes associated to $\alpha$ and $\beta$, choose trivializations $\tau\colon\kx_{\ol \bk}\to B\IG_{m,\ol C}$ and $\tau'\colon\mathscr Y_{\ol \bk}\to B\IG_{m,\ol {C'}}$ as in Lemma~\ref{lem:trivializations}(1), and put
\[
F:=\theta_{\tau'}^{-1}\circ R_\Phi\circ\theta_\tau\colon J\times\hat J\longrightarrow J'\times\hat {J'} .
\]
We want to show that $F$ is symplectic. This can be checked after base change to $\ol\bk$.  By construction, $(\theta_\tau)_{\ol\bk}$ is simply $R_\tau$ and similarly $(\theta_{\tau'})_{\ol\bk}=R_{\tau'}$. In addition, the formation of $R_\Phi$ commutes with extension of the base field (see \cite[7.1--7.6]{olsson2025twisted}). Thus we have

\[
F_{\{\ol\bk\}}=R_{\tau'}^{-1}\circ R_{\Phi_{\ol\bk}}\circ R_\tau.
\]

 Recall (see Example~\ref{ex:Rtau}) that the trivialization $\tau$ induces an equivalence $\tau^*\colon\operatorname{D}^b(\ol C)=\operatorname{D}^b(B\IG_{m,\ol C})^{(1)}\to\operatorname{D}^b(\ol C,\ol\alpha)$, and $R_\tau$  is simply the twisted Rouquier isomorphism associated to $\tau^*$. Similarly $R_{\tau'^*}=R_{\tau'}$. If we consider the equivalence 
 \[
\Phi_0:=(\tau'^*)^{-1}\circ\Phi_{\ol\bk}\circ\tau^*\colon\operatorname{D}^b(\ol C)\to\operatorname{D}^b(\ol {C'}),
\]
since the construction of twisted Rouquier isomorphisms is functorial, we have 
\[R_{\Phi_0}=R_{(\tau^{'*})^{-1}}\circ R_{\Phi_{\ol\bk}}\circ R_{\tau^*}=R_{\tau'}^{-1}\circ R_{\Phi_{\ol\bk}}\circ R_\tau= F_{\{\ol\bk\}}. \]
Now it is immediate to see that $F_{\{\ol\bk\}}$ is symplectic, as it is the (untwisted) Rouquier isomorphism associated to an equivalence between abelian varieties and so symplectic by the work of Orlov \cite{orlov2002derived} (see  Subsection~\ref{subsec:rouquier}).
\end{proof}

\bibliographystyle{alpha}
\bibliography{biblio}
\end{document}